\documentclass{amsart} 
\usepackage{verbatim}
\usepackage{amsmath}  
\usepackage{amsthm}
\usepackage{amssymb}
\usepackage{enumerate}
\usepackage{graphicx}
\newtheorem{theorem}{Theorem}[section]

\newtheorem{lemma}[theorem]{Lemma}
\newtheorem{corollary}[theorem]{Corollary}

\newtheorem{problem}{Problem}

\theoremstyle{definition}

\newtheorem{definition}[theorem]{Definition}

\theoremstyle{remark}
  
\newtheorem{example}[theorem]{Example}{}

\def\spl#1;#2;#3;{{\operatorname{S}}^{#1}(#2,#3)}

\newcommand{\acal}{{\mathcal A}}
\newcommand{\bcal}{{\mathcal B}}
\newcommand{\ccal}{{\mathcal C}}
\newcommand{\dcal}{{\mathcal D}}
\newcommand{\ecal}{{\mathcal E}}
\newcommand{\fcal}{{\mathcal F}}

\newcommand{\pcal}{{\mathcal P}}

\newcommand{\scal}{{\mathcal S}}

\newcommand{\ycal}{{\mathcal Y}}
\newcommand{\zcal}{{\mathcal Z}}

\newcommand{\setm}{\setminus}
\newcommand{\empt}{\emptyset}
\newcommand{\subs}{\subset}

\def\<{\left\langle}
\def\>{\right\rangle}

\def\oo{\omega_1}
\def\br#1;#2;{\bigl[ {#1} \bigr]^ {#2} }

\def\to{\longrightarrow}

\newcommand{\restr}
{\mathop{\hspace{0.01ex}|\hspace*{-0.02ex}{\grave{}}\hspace{0.4ex}}}

\newcommand{\rest}{\restr}

\def\xsp#1;{X(#1)}
\def\xa{\xsp \acal;}
\def\haa#1;#2;{\operatorname{U}_{#1}(#2)}
\def\ha#1;{\haa \acal;#1;}
\def\hn#1;{\haa ;#1;}

\def\ec#1;{\ecal^{#1}}

\def\simple{tree-like}

\def\wf{well-founded}
\def\oot{{\omega}_2}
\def\ooh{{\omega}_3}
\def\csl{chain-closed}
\def\fnu{\acal_0} 
\def\fni{\acal_i} 
\def\feg{\acal_1}

\theoremstyle{plain}

\def\htt{\operatorname{ht}}
\def\lev#1;#2;{\operatorname{I}_{#1}(#2)}
\def\ppl{\stackrel{.}{+}}

\newcommand{\tp}{\operatorname{tp}}

\def\wdt{\operatorname{wd}}

\newcommand{\iclosed}{$\cap$-closed}
\newcommand{\iclosedness}{$\cap$-closedness}

\newcommand{\amfam}[1]{\widehat{#1}} 
\newcommand{\amset}[2]{\amfunname''_{#1}{#2}} 
\newcommand{\amfunname}{\operatorname{k}}

\newcommand{\amseta}[1]{\amset{\acal}{#1}}

\newcommand{\sta}[1]{{#1}^*}

\def\hnu{\sta{\acal_0}} 
\def\heg{\sta{\acal_1}}

\def\hc#1;{\sta{\fcal_{#1}}}
\newcommand{\rk}{\operatorname{rk}}
\newcommand{\rkgood}{$\rk$-good}

\newcommand{\W}{\operatorname{W}}

\def\rklev#1;#2;{\operatorname{R}_{#1}(#2)}

\author[I. Juh\'asz]{Istv\'an Juh\'asz}
\thanks{The first, third and fourth authors were supported by the 
Hungarian National Foundation for Scientific Research grant no. 25745
}   
\address{Alfr{\'e}d R{\'e}nyi Institute of Mathematics}
\email{juhasz@renyi.hu}

\author[S. Shelah]{Saharon Shelah}
\thanks{The second author was supported by the United States 
Israel Binational
Science Foundation, Publication 714}
\address{The Hebrew University of Jerusalem}
\email{shelah@math.huji.ac.il}

\author[L. Soukup]{
Lajos Soukup }
\thanks{
The third author was partially
supported by Grant-in-Aid for JSPS Fellows No.\ 98259 of the Ministry of 
Education, Science, Sports and Culture, Japan}
\address{Alfr{\'e}d R{\'e}nyi Institute of Mathematics }  
\email{soukup@renyi.hu}

\author[Z. Szentmikl\'ossy]{
Zolt\'an Szentmikl\'ossy}
\address{E\"otv\"os University of Budapest} 
\email{zoli@renyi.hu}

\title{A tall space with a small bottom}
\begin{document}
\maketitle

\begin{abstract}
We introduce a general method of constructing locally compact scattered spaces
from certain families of sets and then, with the help of this method, we
prove that if $\kappa^{<\kappa} = \kappa$ then there is such a space of
height $\kappa^+$ with only $\kappa$ many isolated points. This implies
that there is a locally compact scattered space 
of height $\oot$ with $\oo$ isolated points in ZFC, solving an old problem
of the first author.
\end{abstract}

\section{Introduction}\label{sc:intr}
Let us start by recalling that a 
topological space $X$ is called {\em scattered} if every non-empty 
subspace of $X$ has an isolated point and that such a space has a
natural decomposition into levels, the so called Cantor-Bendixson levels. 
The ${\alpha}^{\text{th}}$ 
Cantor-Bendixson level of  $X$ will be 
denoted by $\lev {\alpha};X;$. 
We shall write $\lev
<{\lambda};X;=\bigcup_{{\alpha}<{\lambda}}
\lev {\alpha};X;$.
The height of $X$, $\htt(X)$,
is the least ${\alpha}$ with $\lev {\alpha};X;=\empt$.
The sequence $\<|I_{\alpha}(X)| : \alpha \in \htt(X)\>$ 
is said to be the {\em cardinal sequence} of $X$.
The width of $X$, $\wdt(X)$,
is defined by  
$\wdt(X) = \sup\{\ |\lev {\alpha};X;|:{\alpha}<\htt(X)\}$.

The cardinality of a $T_3$ , in particular of a 
locally compact, scattered  (in short: LCS) space $X$
is at most $2^{\, |\lev ;X;| }$, hence clearly 
$\htt(X)<(2^{|\lev ;X;|})^+$. Therefore  under $CH$ there is no 
LCS space of height ${\omega}_2$ with only countably many isolated points.    
On the other hand, I. Juh\'asz and W.  Weiss, 
\cite[theorem 4]{JW}, proved in ZFC that
for every ${\alpha}<{\omega}_2$ there is a 
LCS space $X$ with $\htt(X)={\alpha}$ and $\wdt(X)={\omega}$. 
The natural question if the existence of
an LCS space of height ${\omega}_2$ with countable width
follows from $\neg CH$
was answered in the negative  by W. Just, who proved,
\cite[theorem 2.13 ]{J1}, that if one adds Cohen reals 
to a model of $CH$ then in the generic extension there are no
LCS spaces of height ${\omega}_2$
and $\wdt(X)={\omega}$. On the other hand, Baumgartner and Shelah
proved it consistent (with $\neg CH$) that such an LCS space exists.

The above mentioned estimate 
$\htt(X)<(2^{|I(X)}|)^+$ is sharp for 
LCS spaces with countably many isolated points : 
it is easy to construct
an LCS space with countable "bottom" and of height ${\alpha}$
for each ${\alpha}< (2^{\omega})^{+}$ 
(see theorem \ref{tm:notch}).
Much less is known about
LCS spaces with $\oo$ isolated points, for example it is a long standing
open problem whether there is, in ZFC, an LCS space of height $\oot$
and width $\oo$. In fact, as was noticed by Juh{\'a}sz in the mid eighties,
even the much simpler question if there is a ZFC example of an LCS space
of height $\oot$ with only $\oo$ isolated points, turned out to be
surprisingly difficult. On the other hand, 
Mart\'{\i}nez, \cite[theorem 1]{Mar}  proved that it is consistent 
that for each ${\alpha}<\ooh$ there is a 
LCS space of height ${\alpha}$ and width $\oo$.
As the main result of the present paper, we shall give an affirmative answer
to the above
question of Juh{\'a}sz: in section \ref{sc:CH} 
we construct, in ZFC, an  
LCS space of height
$\oot$ with $\oo$ isolated points.
Since this space we construct in theorem \ref{tm:zfc} 
has  width $\oot$, the following question remains: 
\begin{problem}
Is there an
LCS space $X$
of height $\oot$ and width $\oo$ in ZFC?
\end{problem}
The methods used in the proof of theorem \ref{tm:zfc}
do not seem to suffice to get LCS  spaces with $\oo$ isolated points 
of arbitrary height $<\ooh$. Thus we have the following problem:
\begin{problem}
Is there, in ZFC, an LCS space with $\oo$ isolated points and of height 
${\alpha}$  for each ${\alpha}<\ooh$? 
\end{problem}

Although one of our main results, theorem 2.19, generalizes to higher
cardinals, it does not seem to suffice to get the analogous result e.g.
for $\ooh$ instead of $\oot$.
If $2^{\omega}\le \oot$ but $2^{\oo}>\oot$ then neither 
theorem \ref{tm:notch} nor theorem \ref{tm:scch} can be applied to get
an  LCS space of height $\ooh$ with only  $\oot$ many isolated points. 
Thus the following version of Juh{\'a}sz' problem remains open:
\begin{problem}
Is there, in ZFC, an   LCS space $X$
of height ${\omega}_3$ having $\oot$ isolated points?
\end{problem}

Let us mention here that the problem of the existence of 
$(\lambda^+,\lambda)$-thin-tall
spaces, i. e. LCS spaces of width $\lambda$ and height $\lambda^+$,
is mentioned in \cite[ Problem 6.4, p.53]{Sh}. However,
it is erroneously stated there that the existence of a 
$(\lambda^+,\lambda)$-thin-tall space follows from $\lambda^{<\lambda} =
\lambda$ or from the existence of a $\lambda^+$-tree.

\section{A space of height $\oot$ and with  $\oo$ isolated points}
\label{sc:CH}

\begin{definition}\label{df:xa}
Given a family  of sets $\acal$ we define the topological space
$\xa=\<\acal,{\tau}_\acal\>$ as follows: 
${\tau}_{\acal}$ is the coarsest topology
in which  the sets $\ha A;=\acal\cap\pcal(A)$ are clopen 
for each $A\in\acal$, 
in other words: $\{\ha A;,\acal\setm \ha A;:A\in\acal\}$
is a subbase for ${\tau}_\acal$.

We shall write $\hn A;$ instead of $\ha A;$ if 
$\acal$ is clear from the context.

\end{definition}

Clearly $\xa$ is a 0-dimensional $T_2$-space. 
A family $\acal$ is  called {\em \wf\ } iff
$\<\acal,\subset\>$ is well-founded.
In this case we can define the rank-function $\rk:\acal\to On$ as usual:
\begin{displaymath}
\rk(A)=\sup\{\rk(B)+1:B\subsetneq A\},
\end{displaymath}
and write $\rklev {\alpha};\acal;=\{A\in\acal:\rk(A)={\alpha}\}$.

The family $\acal$ is said to be {\em \iclosed\ } iff 
$A\cap B\in\acal\cup \{\empt\}$ whenever  $A,B\in\acal$.

It is easy to see that if $\acal$ is \iclosed, 
then a neighbourhood base in $\xsp \acal;$ of $A\in\acal$
 is formed by the sets
\begin{displaymath}
\W(A;B_1,\dots,B_n)=\haa ;A;\setm \bigcup_1^n\haa ;B_i;, 
\end{displaymath}
where $n\in {\omega}$ and $B_i\subsetneq A$ for $i=1,\dots,n$. ( For $n=0$
we have $\W(A)=\haa ;A;$.)

The following simple result enables us to obtain LCS spaces from certain
families of sets. Let us point out, however, 
that not every LCS space is obtainable
in this manner, but we do not dwell upon this because we will not need it.  

\begin{lemma}\label{lm:scat}
Assume that $\acal$ is both \iclosed\ and well-founded. 
Then $\xa$ is an LCS space. 
\end{lemma}

\begin{proof}

Given  a non-empty subset $\ycal$ of $\acal$ let $A$ be a $\subset$-minimal
element of $\ycal$. Then $\hn A;\cap\ycal=\{A\}$, i.e. $A$ is isolated in
$\ycal$. Thus $\xa$ is scattered.

Next we prove that every $\hn A;$ is compact by well-founded induction on 
$\<\acal,\subset\>$. Assume that  $\hn B;$ is compact
for each $B\subsetneq A$. By Alexander's subbase lemma it is enough
to prove that any cover of $\hn A;$ with subbase elements  contains a
finite subcover.
So let $\bcal,\ccal\subs \acal$ be such that 
\begin{equation}\notag
\hn A;\subs \bigcup_{B\in\bcal}\hn B;\cup\bigcup_{C\in\ccal}(\acal\setm \hn C;).
\end{equation}                                           
If $A\in\hn B;$ for some $B\in\bcal$ 
then $A\subs B$ and so $\hn A;\subs \hn B;$.

Hence we can assume that $A\in\acal\setm \hn C;$, i.e.
$A\not\subs C$ for some $C\in\ccal$. 
If $A\cap C=\empt$ then $\hn A;\setm \{\empt\}\subs \acal\setm \hn C;$,
and we are clearly done. 
So we can assume that $A\cap C\ne \empt$, and 
consequently $A\cap C\in\acal$. Then  
$\hn A;\setm (\acal\setm \hn C;)=\hn A;\cap \hn C;=\hn A\cap C;$.
Since $A\cap C\ne A$ the set $\hn A\cap C;$ is compact 
by the induction hypothesis,
hence $\hn A\cap C;$ is covered by a finite subfamily $\fcal$ of 
$\{\hn B;:B\in\bcal\}\cup\{\acal\setm \hn D;:D\in\ccal\}$.
Therefore $\fcal\cup\{\acal\setm \hn C;\}$ is a finite cover of $\hn A;$. 
Consequently $\hn A;$ is compact. 
\end{proof}

To simplify notation, if $\xsp\acal;$ is scattered then we  write
$\lev {\alpha};\acal;=\lev {\alpha};{\xsp\acal;};$.

Clearly each minimal element of $A\in\acal$ is isolated in $\xsp \acal;$;
more generally we have ${\alpha}\le \rk(A)$ if $A\in \lev {\alpha};\acal;$,
as is shown by an easy induction on $\rk(A)$.

\begin{example}\label{ex:kappa} 
Assume that $\<T,\prec\>$ is a well-ordering,
$\tp \<T,\prec\>={\alpha}$, and let
$\acal$ be the family of all initial segments of $\<T,\prec\>$, i. e.
$\acal=\{T\}\cup \{T_x:x\in T\}$, where $T_x=\{t\in T:t\prec x\}$.
Then $\acal$ is \wf, \iclosed\ and it is easy to see that
$\xsp \acal;\cong {\alpha}+1$, i.e. the space $\xsp \acal;$ is homeomorphic
to the space of ordinals up to and including $\alpha$.

\end{example}

Example \ref{ex:kappa} above shows that,
in general,
$\rklev {\alpha};\acal;$ and $\lev {\alpha};\acal;$ may differ even
for ${\alpha}=0$. Indeed,
if $x$ is the successor of $y$ in $\<T,\prec\>$ then $T_x$ is isolated 
in $\xsp\acal;$
because
$\{T_x\}=W(T_x;T_y)=\ha T_x;\setm \ha T_y;$ is open, but 
$\rk(x)=\tp(T_x)>0$.
However, for a wide class of families, the two kinds of levels do agree. 
Let us call a \wf\ family $\acal$ {\em \rkgood}
iff the following condition is satisfied:
\begin{displaymath}
\forall A\in\acal\  \forall {\alpha}<\rk(A)\ 
|\{A'\in \acal:A'\subs A\land \rk(A')={\alpha}\}|\ge {\omega}. 
\end{displaymath}

Then we have the following result.

\begin{lemma}
If $\acal$ is \wf, \iclosed\ and \rkgood\ then
$\lev {\alpha};\acal;=\rklev {\alpha};\acal;$ for each ${\alpha}$.
\end{lemma}

\begin{proof}
We prove this by induction on ${\alpha}$. 
Assume that $\lev {\xi};\acal;=\rklev {\xi};\acal;$ for all ${\xi}<{\alpha}$.

If $A\in \rklev {\alpha};\acal;$ then
$\haa ;A;\setm \{A\}\subs \bigcup_{{\xi}<{\alpha}}\rklev {\xi};\acal;=
\bigcup_{{\xi}<{\alpha}}\lev{\xi};\acal;$ and so 
$A$ is an isolated point of  $\acal\setm \bigcup_{{\xi}<{\alpha}}\lev{\xi};\acal;$, i.e.
$A\in \lev {\alpha};\acal;$. Thus we have $\rklev
{\alpha};\acal;\subs \lev {\alpha};\acal;$.

Now assume that $A\in\lev{\alpha};\acal;\setm \rklev {\alpha};\acal;$. 
Then by our above remark ${\alpha}<\rk(A)$, moreover there are $B_1,\dots
B_n\in\ha A;\setm \{A\}$ such that 
\begin{equation}\tag{$\star$}\label{eq:iso}
\W(A;B_1,\dots,B_n)\setm\{A\}\subs \lev <{\alpha};\acal;
=\rklev <{\alpha};\acal;.
\end{equation}
Let ${\eta}=\max \{{\alpha},\max_{i=1,\dots,n}\rk B_i\}$.
Then  ${\eta}<\rk (A)$, moreover  
we have $\ha A;\cap \rklev {\eta};\acal;\subs \{B_1,\dots B_n\}$
by (\ref{eq:iso}), contradicting  
$|\ha A;\cap \rklev {\eta};\acal;|\ge{\omega}$.
Note that this argument is valid for $n=0$ as well. Indeed, in this case
we have $\haa \acal;A;\setm\{A\}\subs \lev <{\alpha};\acal;$,  
moreover ${\eta}={\alpha}$.
Thus we have concluded
that $\lev {\alpha};\acal;=\rklev {\alpha};\acal;$.
\end{proof}

\begin{example}\label{exa}
For a fixed  cardinal ${\kappa}$ and any ordinal ${\gamma}<{\kappa}^+$
we define the family 
$\ec {\gamma};\subs\pcal({\kappa}^{\gamma})$
as follows:
\begin{displaymath}
\ec {\gamma};=\Bigl\{ \bigl[ 
{\kappa}^{1+\alpha}\cdot{\xi},{\kappa}^{1+\alpha}\cdot({\xi}+1)
 \bigl) \ : \  {\alpha}\le {\gamma}, \  {\kappa}^{1+\alpha}\cdot{\xi}< 
{\kappa}^{\gamma}
\Bigr\}. 
\end{displaymath}

Of course, throughout this definition exponentiation means ordinal exponentiation.

$\ec {\gamma};$ is clearly \wf, \iclosed, moreover
$\rk\Bigl(\bigl[ 
{\kappa}^{1+\alpha}\cdot{\xi},{\kappa}^{1+\alpha}\cdot({\xi}+1)
 \bigl)\Bigr)={\alpha}$, hence $\ec {\gamma};$ is also \rkgood. Consequently 
$\xsp {\ec {\gamma};};$ is an LCS
space of height ${\gamma}+1$ in which  
the ${\alpha}^{\text{th}}$ level is
$\Bigl\{\big[ 
{\kappa}^{1+\alpha}\cdot{\xi},{\kappa}^{1+\alpha}\cdot({\xi}+1
)\big) \ :   {\kappa}^{1+\alpha}\cdot{\xi}< {\kappa}^{\gamma} \Bigr\}$,
i. e. all levels except the top one are of size $\kappa$.
\end{example}

To get an LCS space of height ${\kappa}^+$ with ``{\em few}'' isolated
points, 
our plan is to {\em amalgamate}
the spaces $\{\xsp{\ec {\gamma};};:{\gamma}<{\kappa}^+\}$ 
into one LCS space $X$ in such a way that 
$|\lev 0;X;|\le {\kappa}^{<{\kappa}}$.
The following definition describes a situation in which such an
amalgamation can be done.

\begin{definition}
A system of families $\{\acal_i:i\in I\}$ is called {\em coherent}
 iff $A\cap B\in \acal_i\cup\{\empt\}$ whenever
$\{i,j\}\in \br I;2;$, $A\in \acal_i$ and $B\in \acal_j$.
\end{definition}

To simplify notation, we introduce the following convention.
Whenever the system of families $\{\acal_i:i\in I\}$ is given,
we will write $\haa i;A;$ for $\haa \acal_i;A;$, and
${\tau}_i$ for ${\tau}_{\acal_i}$. If the family $\acal$ is defined
then we will write $\haa ;A;$ for $\haa \acal;A;$,
and ${\tau}$ for ${\tau}_{\acal}$.

\begin{lemma}\label{lm:coherent}
Assume that $\{\acal_i:i\in I\}$ is a coherent system of \wf, \iclosed\ 
families and $\acal=\cup\{\acal_i:i\in I\}$. Then
for each $i\in I$ and $A\in\acal_i$ we have
$\haa ;A;=\haa i;A;$, 
$\acal$ is also \wf\and  \iclosed, moreover 
${\tau}_{i}\rest \haa ;A;={\tau}\rest
\haa ;A;$.
Consequently each $\xsp \acal_i;$ is an open subspace
of $\xsp \acal;$  
and thus $\{\xsp \acal_i;:i\in I\}$ forms  an open cover of
$\xsp \acal;$.
\end{lemma}

\begin{proof}
Let $A\in\acal_i$.
Then it is clear from coherence that
\begin{equation}\notag
\haa i;A;\subseteq\haa ;A;=
\bigcup_{j\in I}\{B:B\in \acal_j\land B\subs A\}
\subseteq
\haa i;A;,
\end{equation}
hence $\haa i;A;=\haa ;A;$.

Next let $B\in\acal_j$.
If $A\cap B=\empt$ then $\haa ;A;\cap \haa ;B;\subset\{\empt\}.$
Now assume that $A\cap B\ne \empt$.
Then, again by coherence, $A\cap B\in \acal_i$ and we have
\begin{multline}\notag
\haa ;A;\cap \haa ;B;=\haa i;A;\cap\haa ;B;
=\{C\in\acal_i:C\subs A\land C\subs B\}\\
=\{C\in\acal_i:C\subs A\cap B\}
=\haa i;A\cap B;.
\end{multline}

In both cases $\haa ;A;\cap \haa ;B;$ is ${\tau}_i$-open.
Similarly we can see that $\haa ;A;\setm \haa ;B;=\haa i;A;\setm \haa ;B;$
is ${\tau}_i$-open,
hence the topologies
${\tau}_i\rest \haa ;A;$ and 
${\tau}_{}\rest\haa ;A;$ coincide.

To show that $\acal$ is {\wf},
assume  that $\{A_n:n\in {\omega}\}\subs \acal$ and 
$A_0\supseteq A_1\supseteq \dots$.
If $A_0\in \acal_i$ then $\{A_n:n\in {\omega}\}\subs \acal_i$
because $\haa ;A_0;=\haa i;A_0;$. Thus
there is $n\in {\omega}$ with $A_m=A_n$ 
for each $m\ge n$ because 
$\acal_i$ is \wf. Finally, that $\acal$ is \iclosed\ is an easy consequence of
coherence and the \iclosedness of the families $\acal_i$.
\end{proof}

Given a system of families $\{\acal_i:i\in I\}$
we would like to construct a coherent system of families 
$\{\amfam{\acal_i}:i\in I\}$ such that $\acal_i$ and $\amfam{\acal_i}$
are isomorphic for all $i \in I$. 
A sufficient condition for when this can be done will be given in lemma
\ref{lm:am2} below. 

First, however, we need a definition.
While reading it, one should remember that an ordinal
is identified with the family of its proper initial segments.

\begin{definition}\label{df:am1}
Given a limit ordinal  ${\rho}$ and a family $\acal$ with 
${\rho}\subs\acal\subs \pcal({\rho})$, let us
define the family $\amfam {\acal}$ as follows.
Consider first the function $\amfunname_{\acal}$ on ${\rho}$ determined
by the formula
$\amfunname_{\acal}({\eta})=\haa
\acal;{\eta}+1;$ for $\eta \in \rho$ and  
put
\begin{displaymath}
\amfam{\acal}=\{\amseta A:A\in\acal\}.
\end{displaymath}

Since  ${\rho}\subs\acal$, for each ${\eta}\in {\rho}$ we clearly have 
$\cup \ha {\eta};={\eta}$ and so
$\amfunname_{\acal}({\eta})=$
$\ha {\eta}+1;\ne \ha {\xi}+1;=\amfunname_{\acal}({\xi})$
whenever $\{{\eta},{\xi}\}\in \br {\rho};2;$.
Consequently, 
$\amfunname_{\acal}$ is a bijection that yields an isomorphism between
$\acal$ and $\amfam{\acal}$ (and so the spaces
$\xsp \acal;$ and $\xsp{\amfam{\acal}};$ are homeomorphic). 

If the system of families $\{\acal_i:i\in I\}$ is given, then
we write $\amfunname_i$ for $\amfunname_{\acal_i}$ for each $i \in I$.
\end{definition}

If $\acal \subs \pcal({\rho})$ and ${\xi}\le {\rho}$ then we let
\begin{displaymath}
\acal\rest {\xi}=\{A\cap {\xi}:A\in \acal\}.
\end{displaymath}

For $\fnu\ne\feg\subs \pcal({\rho})$ we let
\begin{displaymath}
\Delta(\fnu,\feg)=\min\{{\delta}:\fnu\rest{\delta}\ne 
\feg \rest {\delta}\}.
\end{displaymath}

Clearly we always have $\Delta(\fnu,\feg)\le {\rho}$.
If, in addition,  
${\rho}+1\subs \fnu\cap\feg$, moreover both $\fnu$ and $\feg$ are 
\iclosed\, then we also have
\begin{displaymath}
\Delta(\fnu,\feg)=\min\{{\delta}: \haa 0;{\delta};
\ne \haa 1;{\delta};\}, 
\end{displaymath}
because then $\acal_i\rest {\delta}=\haa i;{\delta};$
whenever $i \in 2$ and ${\delta}\le {\rho}$.

\begin{lemma}\label{lm:am2}
Assume that  ${\kappa}$ is a cardinal, 
$\{\acal_i:i\in I\}\subs \pcal\pcal({\kappa})$
are \iclosed\ families,
${\kappa}+1\subs \acal_i$ for each $i\in I$, and $\Delta(\acal_i,\acal_j)$
is a successor ordinal whenever $\{i,j\}\in \br I;2;$. Then the system
$\{\amfam{\acal_i}:i\in I\}$ is coherent.
\end{lemma}

\begin{proof}

Let $A\in\acal_i$ and $B\in\acal_j$,
where ${\Delta}(\acal_i,\acal_j)=
{\rho}+1$. 
 Then  $B\cap {\rho}\in \haa j;{\rho};=\haa
i;{\rho};$ by the choice of ${\rho}$
and so $A\cap B\cap {\rho}\in {\acal_i}\cup \{\empt\}$ because
$\acal_i$ is \iclosed. Consequently  we have
\begin{multline}\notag
\amset {i}A\cap \amset {j}B=\\
\bigl\{\haa i;{\eta}+1;:
{\eta}\in A\cap B\ \land\ \haa i;{\eta}+1;=
\haa j;{\eta}+1; \bigr\}=
\\\bigl\{\haa i;{\eta}+1;:
{\eta}\in A\cap B\ \land\ {\eta}<{\rho}\bigr\}=
\amset {i}({A\cap B\cap {\rho}})\in \amfam{\acal_i}
\cup\{\empt\} , 
\end{multline}
as required by the definition of coherence.
\end{proof}

More is needed still if we want the "amalgamated" family to provide us
a space with a small base, i.e. having not too many isolated points.
This will be made clear by the following lemma. 

\begin{lemma}\label{lm:amalg}
Let ${\kappa}$ be a cardinal and    
$\{\acal_i:i\in I\}\subs \pcal\pcal({\kappa})$ be a system of 
families
such that 
\begin{enumerate}[\rm (i)]
\item\label{ai} ${\kappa}\ppl 1\subs \acal_i$ 
and $\acal_i$ is {\wf} and \iclosed\ 
for each $i\in I$,
\item\label{succ}  
${\Delta}(\acal_i,\acal_j)
$ is a successor ordinal for each  $\{i,j\}\in \br I;2;$.
\end{enumerate}
Then 
\begin{enumerate}[\rm (a)]
\item \label{ama} 
the system $\{\amfam{\acal_i}:i\in I\}$ is coherent 
and thus 
$\acal=\bigcup\{\amfam{\acal_i}:i\in I\}$
 is {\wf}, \iclosed\ 
and
$\xsp \acal;$ is covered by its open subspaces 
$\{\xsp {\amfam{\acal_i}};:i\in I\}$.
\end{enumerate}
If, in addition, we also have
\begin{enumerate}[\rm (i)]\addtocounter{enumi}{2}
\item \label{aiat} 
$\lev 0;\acal_i;\subs \br {\kappa};<{\kappa};$ for each $i\in I$, 
\end{enumerate}
and
\begin{enumerate}[\rm (i)]\addtocounter{enumi}{3}
\item \label{ueta}
$|\haa{i};{\eta};|<{\kappa}$ 
for each $i\in I$ and ${\eta}\in {\kappa}$,
\end{enumerate}
then 
\begin{enumerate}[\rm (a)]\addtocounter{enumi}{1}
\item \label{amalat}
$\lev 0;\acal;|\subs 
\bigg[\Big[{[{\kappa}]^{<{\kappa}}}\Big]^{<{\kappa}}\bigg]^{<{\kappa}}$.
\end{enumerate}
\end{lemma}

\begin{proof}[Proof of lemma \ref{lm:amalg}]

The system $\{\amfam{\acal_i}:i\in I\}$ is coherent 
by lemma \ref{lm:am2}, thus (\ref{ama}) holds by lemma \ref{lm:coherent}.  
   
Consequently we have
\begin{displaymath}
\lev 0;{\acal};= \bigcup 
\{\lev 0;{\amfam {\acal_i}};:i\in I\}.
\end{displaymath}

Now if $A\in\lev 0;\acal_i;$ and ${\eta}\in {\kappa}$ then  
$|A|<{\kappa}$ by (\ref{aiat})
and $\haa i;{\eta};\in \Big[\br {\kappa};<{\kappa};\Big]^{<{\kappa}}$
by (\ref{ueta}), hence 
\begin{displaymath}
\amset{i}{A}=\{\haa i;{\eta}+1;:{\eta}\in A\}\in \bigg[\Big[{\br {\kappa};<{\kappa};}\Big]^{<{\kappa}}
\bigg]^{<{\kappa}}.
\end{displaymath}

This, by $\lev 0;{\amfam{\acal_i}};=\{\amset{i}{A}:A\in \lev 0;\acal_i;\}$,
proves (\ref{amalat}).
\end{proof}

Before we could apply this result to the families $\ec {\gamma};$,
however, we need some further preparation.

\begin{definition}\label{df:simple}
A family $\acal$ is called {\em  \simple}\ iff 
$A\cap A'\ne\empt$ implies that  $A\subs A'$ or $A'\subs A$,
whenever  $A,A'\in\acal$.
\end{definition}

\begin{definition}\label{df:ccl}
A family $\acal$ is called {\em\csl} if for each non-empty
$\bcal\subs \acal$ if $\bcal$ is ordered by $\subs $ (i.e. if
$\bcal$ is a chain) then
$\cup\bcal\in\acal$.
\end{definition}

It is easy to see that the families $\ec {\gamma};$ given in example 
\ref{exa} are both \simple\ and \csl.
Also, \simple\ families are clearly \iclosed.

\begin{lemma}\label{lm:rest}
If ${\delta}$ is an ordinal and $\acal\subs \pcal({\delta})$ is \simple, {\wf} and
\csl\ then so is $\acal\rest {\xi}$ for each ${\xi\def\spl#1;#2;#3;{{\mathfrak S}^{#1}(#2,#3)}
}\le {\delta}$. 
\end{lemma}

\begin{proof}[Proof of lemma \ref{lm:rest}]
It is obvious that  $\acal\rest {\xi}$ is \simple. 
To show that $\acal\rest {\xi}$ is \csl, let
 $\empt\ne \bcal \subs \acal\rest {\xi}$  be ordered by $\subs$ . 
If $\bcal=\{\empt\}$ then $\cup\bcal=\empt\in\acal\rest {\xi}$.
If, however, $\bcal\ne \{\empt\}$ then 
put $\tilde{\bcal}=\{A\in\acal:A\cap {\xi}\in \bcal\setm \{\empt\}\}$.  
Since $\acal$ is \simple, $\tilde{\bcal} $ is also ordered by $\subs$
and clearly $\tilde{\bcal}\ne\empt$.
So $\cup \tilde{\bcal} \in \acal$ and 
$\cup\bcal=\cup \tilde{\bcal} \cap {\xi}\in \acal\rest {\xi}$, 
which was to be shown.

To show that $\acal\rest {\xi}$ is {\wf}
assume  that $A_0\cap {\xi}\supseteq A_1\cap {\xi}\supseteq \dots$, 
where each $A_n\in\acal$.
If $A_n\cap {\xi}=\empt$ for some $n$, then we are done.
Otherwise for each $n\in {\omega}$ we have 
$A_n\cap {\xi}=(\cap_{m\le n}A_m)\cap {\xi}\ne\empt$, hence
as $\acal$ is \iclosed\ we can assume
that $A_0\supseteq A_1\supseteq\dots$. Since $\acal$ is \wf, there is
$n$ such that $A_m=A_n$,  and so
$A_m\cap {\xi}=A_n\cap {\xi}$ as well, for each $m\ge n$. 
\end{proof}

\begin{definition}\label{df:ab}
Given a family $\acal\subs \pcal({\delta})$
and ${\alpha},{\beta}\in {\delta}$ let us put
\begin{equation}\notag
\spl \acal;{\alpha};{\beta};=\cup\{A\in\acal:{\alpha}\in
A\text{ and }{\beta}\notin A\}. 
\end{equation}
\end{definition}

\begin{lemma}\label{lm:char} 
Assume that ${\delta}$ is an infinite ordinal and 
$\acal\subs \pcal({\delta})$ is a \simple, \wf\ and
\csl\ family with ${\delta}\in \acal$. 
Then  
\begin{equation}\notag
\acal\setm \{\empt\}=\{{\delta}\}\cup
\{\spl \acal;{\alpha};{\beta};:{\alpha},{\beta}\in {\delta}\}\setm \{\empt\}.
\end{equation}
Consequently, $|\acal|\le |{\delta}|$.
\end{lemma}

\begin{proof}[Proof of lemma \ref{lm:char}]
Given ${\alpha},{\beta}\in {\delta}$, the family 
$\scal=\{A\in\acal:{\alpha}\in A, {\beta}\notin A\}$ is ordered by 
$\subs$ because $\acal$ is \simple. Thus
either $\scal=\empt$ and so $\spl \acal;{\alpha};{\beta};=\cup\scal=\empt$,
or 
 if $\scal\ne\empt$ then $\spl
\acal;{\alpha};{\beta};=\cup\scal\in\acal$, for $\acal$ is \csl.

Assume now that $A\in \acal\setm \{{\empt,\delta}\}$ and 
let $\dcal=\{D\in\acal:A\subsetneq D\}$. Clearly ${\delta}\in \dcal$.
Since $\acal$ is \simple, $\dcal$ is  ordered
by $\subs$ , so it has a $\subs$-least element, say $D$, because 
$\<\acal,\subs\>$ is also \wf.
Pick ${\beta}\in D\setm A$ and let ${\alpha}\in A$.
We claim that  $A=\spl \acal;{\alpha};{\beta};$.
Clearly $A\subs\spl \acal;{\alpha};{\beta};$ because ${\alpha}\in A$
and ${\beta}\notin A$. On the other hand, if $A'\in\acal$, 
${\alpha}\in A'$ and 
${\beta}\notin A'$ then either $A'\subs A$ or $A\subset A'$
because $\acal$ is \simple. But ${\beta}\notin A'$ implies that 
$A'\notin\dcal$, i.e. $A\subsetneq A'$ can not hold.
Thus $A'\subs A$ and so $\spl \acal;{\alpha};{\beta};= A$ is proved.
\end{proof}

\begin{definition}\label{df:sta}
If ${\rho}$ is an ordinal and $\acal \subs \pcal({\rho})$  let us put
\begin{displaymath}
\sta{\acal}=
\{A\cap {\xi}:A\in\acal\land {\xi}\le {\rho}\}=
\acal\cup\{A\cap {\xi}:A\in\acal\land {\xi}<{\rho}\}.
\end{displaymath}
\end{definition}

\begin{lemma}\label{lm:dpres}
If ${\rho}$ is an ordinal and $\fnu,\feg\subs \pcal({\rho})$ are \csl,
\iclosed\ and  \wf\ families 
such that $\hnu\ne\heg$
then  $\Delta(\hnu,\heg)$ is a successor ordinal.
\end{lemma}

\begin{proof}
Assume that ${\delta}$ is a limit ordinal and
$\sta \fnu\rest {\gamma}=\sta\feg\rest {\gamma}$ for  all ${\gamma}<{\delta}$.
We want to show that $\sta \fnu\rest {\delta}=\sta\feg\rest {\delta}$.
Since $\sta\fni\rest {\delta}=\bigcup_{{\xi}\le {\delta}}\fni\rest {\xi}$
and $\bigcup_{{\xi}<{\delta}}\fni\rest {\xi}=\bigcup_{{\xi}<
{\delta}}\sta\fni\rest {\xi}$, moreover $\empt\in \sta \fnu\cap\sta\feg$,
it is enough to show that 
$(\fnu\rest {\delta})\setm \{\empt\}=(\feg\rest {\delta})\setm \{\empt\}$. 

So assume that 
$A\in \fnu$ 
with $A\cap {\delta}\ne \empt$
and verify  that then $A\cap {\delta}\in\feg\rest {\delta}$.

Fix ${\zeta}\in A\cap {\delta}$. 
For each ${\gamma}$ with 
${\zeta}<{\gamma}<{\delta}$ let $B_{\gamma}$ be the $\subset$-minimal
element of $\feg$ with $B_{\gamma}\cap {\gamma}=A\cap {\gamma}\ne \empt$.
Then $\{B_{\gamma}:{\gamma}<{\delta}\}$ is a chain because
$B_{\gamma}\subs B_{\gamma'}$ for ${\gamma}<{\gamma}'$
by the minimality of $B_{\gamma}$ and because $\feg$ is \iclosed. Thus
$B=\cup\{B_{\gamma}:{\zeta}<{\gamma}<{\delta}\}\in\feg$
and clearly $A\cap {\delta}=B\cap {\delta}$.
\end{proof}

The last result shows us that the operation * is useful because its
application yields us families that satisfy condition (ii) of lemma
2.10. On the other hand, the following result tells us that the LCS
spaces associated with the families modified by * 
do not differ significantly from the spaces given by
the original families, moreover they also satisfy condition (iii)
of lemma 2.10.    

\begin{lemma}\label{lm:kappa}
Let ${\kappa}$ be a cardinal and $\acal\subs \br {\kappa};{\kappa};$ be {\wf} and \iclosed. 
Then so is $\sta\acal$, moreover
\begin{enumerate}[$(a)$]
\item \label{cl}
$\xsp \acal;$ is a closed subspace of $\xsp {\sta\acal};$,
\item \label{kappa} $\lev 0; \acal;\subseteq\lev {\kappa};{\sta\acal};$.
\item \label{htk} $\htt({{\sta\acal}})\ge {\kappa}\ppl\htt( \acal)$,
\item\label{le} $\lev 0;{\sta\acal};\subs
\br{\kappa};<{\kappa};$.
\end{enumerate}
\end{lemma}

\begin{proof}[Proof of lemma \ref{lm:kappa}]
We shall write $\haa ;A;$ for $\haa \acal;A;$,
and $\haa *;A;$ for $\haa \sta{\acal};A;$.

First observe that because
\begin{displaymath}
\haa{*};{A};\cap\acal=
\left\{ 
\begin{array}{ll}
\haa ;A;&\text{if $A\in\acal$,}\\
\emptyset&\text{if $A\in\sta\acal\setm \acal$,}
\end{array}
\right.
\end{displaymath}
$\xsp\acal;$ is a closed subspace of $\xsp{\sta\acal};$,
 hence (\ref{cl}) holds.

Now let $A\in \lev 0;{\acal};$. Then there are 
$B_1,\dots, B_n\in\haa ;A;\setm \{A\}$
such that 
\begin{displaymath}
\{A\}=\W(A;B_1,\dots, B_n)=
\haa ;A;\setm \bigcup_{i=1}^n\haa ;B_i;.
\end{displaymath}
Since here $B_i\subsetneq A$ and $|A| = \kappa$, we can 
fix ${\eta}\in A$ such that $(A\cap {\eta})\not\subset B_i$
for every $i=1,\dots, n$.

Now consider the basic neighbourhood
\begin{displaymath}
\zcal =\W_*(A;A\cap {\eta},B_1,\dots, B_n)=
\haa *;A;\setm \haa *;A\cap {\eta};\setm\bigcup_{i=1}^n\haa *;B_i;
\end{displaymath}
of $A$ in $\xsp{\sta{\acal}};$. We claim that 
$\zcal =\{A\cap {\xi}:{\eta}< {\xi}\le {\kappa}\}$.
The inclusion $\supset$ is clear from the choice of ${\eta}$.
On the other hand, if $C\cap {\xi}\in \zcal$ with $C\in\acal$ and ${\xi}\le
{\kappa}$, 
then $C\cap {\xi}\subs A$ hence  
$C\cap {\xi}=A\cap C\cap {\xi}$, so as $\acal$ is \iclosed\ 
we can assume that $C\subseteq A$.
If we had $C\ne A$ then 
$\{A\}=\W(A;B_1,\dots, B_n)$ would imply 
$C\subs B_i$ for some $i$, hence 
$C\cap {\xi}\in \haa *;B_i;$ and so $C\cap {\xi}\notin \zcal$,
a contradiction,
thus we must have $C=A$. 
Moreover, since $\haa *;A\cap {\eta};\supset \{A\cap {\nu}:{\nu}\le {\eta}\}$,
we must also have ${\xi}> {\eta}$.

By example \ref{ex:kappa} we have
$\xsp \zcal ;\cong {\kappa}+1$.
Moreover, the topologies ${\tau}_{\zcal }$ and ${\tau}_{\sta{\acal}}\rest \zcal $
coincide because the above argument also shows that 
for each $C\in\acal$ and ${\zeta}\le {\kappa}$ we have
\begin{displaymath}
\haa{*};{C\cap {\zeta}};\cap\zcal=
\haa{*};{A\cap C\cap {\zeta}};\cap\zcal=
\left\{ 
\begin{array}{ll}
\haa \zcal;A\cap {\zeta};&\text{if $A\subs C$ and ${\zeta}>{\eta}$;}\\
\emptyset&\text{otherwise.}
\end{array}
\right.
\end{displaymath}

Hence $\xsp \zcal;\cong {\kappa}\ppl 1$ is a clopen subspace of $\xsp {\sta\acal};$ and so 
$\{A\}=\lev {\kappa};\zcal;=\lev {\kappa};{\sta\acal};\cap\zcal$,
what proves (\ref{kappa}).

(\ref{htk}) follows immediately from (\ref{cl}) and (\ref{kappa}).

Finally, $\lev 0;{\sta\acal};
\subs\lev <{\kappa};{\sta{\acal}};
\subs (\sta\acal\setm \acal)
\subs \br {\kappa};<{\kappa};$,  as follows immediately from 
(\ref{kappa}), proving (\ref{le}).
\end{proof}

Now we are ready to collect the fruits of all the preparatory work.

\begin{theorem}\label{tm:scch}
If ${\kappa}^{<{\kappa}}={\kappa}$ 
then there is an LCS space $X$
of height ${\kappa}^+$ with $|\lev 0;X;|={\kappa}$.
\end{theorem}

\begin{proof}[Proof of theorem \ref{tm:scch}]
For each ${\gamma}<{\kappa}^+$ consider the \wf, \iclosed, \rkgood\ family 
$\ec {\gamma};$ constructed in example \ref{exa}:
\begin{displaymath}
\ec {\gamma};=\Bigl\{ \bigl[ 
{\kappa}^{1+\alpha}\cdot{\xi},{\kappa}^{1+\alpha}\cdot({\xi}+1)
 \bigl) \ : \  {\alpha}\le {\gamma}, \  {\kappa}^{1+\alpha}\cdot{\xi}< 
{\kappa}^{\gamma}
\Bigr\}. 
\end{displaymath}
Fix a bijection $f_{\gamma}:{\kappa}^{\gamma}\to {\kappa}$, and let  
$\fcal_{\gamma}=\{f_{\gamma}{}''E:E\in \ec {\gamma};\}$, i.e.
$\fcal_{\gamma}$ is simply an isomorphic copy of $\ec {\gamma};$ on
the underlying set $\kappa$. As
$\ec {\gamma};$ is also \csl\ and \simple, hence so is $\fcal_{\gamma}$.

We shall now show that the *-modified families 
$\{\sta{\fcal_{\gamma}}:{\gamma}<{\kappa}^+\}$
satisfy  conditions (\ref{ai})-(\ref{ueta})
of lemma \ref{lm:amalg}.
Since ${\kappa}\in \fcal_{\gamma}$ it follows that 
${\kappa}\ppl 1\subs \hc {\gamma};$ and so 
(\ref{ai}) is  true.
For $\{{\gamma},{\delta}\}\in \br {\kappa}^+;2;$,
the height of $\xsp {\ec {\gamma};};$ is ${\gamma}+1$
and the height of $\xsp {\ec {\delta};};$ is ${\delta}+1$,
hence $\ec {\gamma};$ and $\ec {\delta};$ are not isomorphic. Thus
$\fcal_{\gamma}\ne \fcal_{\delta}$ and so
$\sta{\fcal_{\gamma}}\ne \sta{\fcal_{\delta}}$
because $\fcal_{\gamma}=\sta{\fcal_{\gamma}}\cap \br {\kappa};{\kappa};$ 
and $ \fcal_{\delta}=\sta{\fcal_{\delta}}\cap \br {\kappa};{\kappa};$.
Hence 
$\Delta(\hc {\gamma};,\hc {\delta};)$ is a successor ordinal 
by lemma \ref{lm:dpres}, i.e. (ii) is satisfied.

(\ref{aiat}) holds by \ref{lm:kappa}.(\ref{le}.)

To show (\ref{ueta}), let us fix ${\xi}<{\kappa}$. Then 
$\haa {\sta{\fcal_{\gamma}}};{\xi};=\hc {\gamma};\rest {\xi}=\cup\{\fcal_{\gamma}\rest {\zeta}:{\zeta}\le {\xi}\}$
 where $|\fcal_{\gamma}\rest {\zeta}|\le |{\zeta}|^2$ for all ${\zeta}\le {\xi}$
by lemmas \ref{lm:rest} and \ref{lm:char}, consequently
$|\hc {\gamma};\rest {\xi}|\le |{\xi}|^3<{\kappa}$.
 
Thus we may apply lemma \ref{lm:amalg} to the family  
$\fcal=\cup\{\amfam{{\hc {\gamma};}}:{\gamma}<{\kappa}^+\}$
and conclude that the
space  $X=\xsp{\fcal};$ 
is LCS,
$|\lev 0;X;|\le \big(({\kappa}^{<{\kappa}})^{<{\kappa}}\big)^{<{\kappa}}
={\kappa}$,
moreover since for every ${\gamma}\in {\kappa}^+$
the space $\xsp {\hc {\gamma};} ;$ is an open subspace
of $X$, we have  
$\htt (X)\ge \htt (\xsp {\hc {\gamma};};)>
{\gamma}$,
consequently $\htt(X)\ge {\kappa}^+$.
\end{proof} 

In particular, if $2^{\omega}=\oo$ then the above result yields 
an LCS space $X$
with $\htt(X)={\omega}_2$ and $|\lev 0;X;|=\oo$. That such a space also exist
under $\neg {}$CH, hence in ZFC, follows from the following result.

\begin{theorem}\label{tm:notch}
For each ${\alpha}<(2^{\omega})^+$ there is a 
locally compact, scattered space $X_{\alpha}$ with 
$|X_{\alpha}|\le |{\alpha}|+{\omega}$, 
$\htt(X_{\alpha})={\alpha}$ and $|\lev 0;X_{\alpha};|={\omega}$.
\end{theorem}

\begin{proof}
We do induction on ${\alpha}$.
If ${\alpha}={\beta}+1$ then we 
let $X_{\alpha}$ be the 1-point compactification
of the disjoint topological sum of countably many copies of $X_{\beta}$.

If ${\alpha}$ is limit then we first fix
an almost disjoint family
$\{A_{\beta}:{\beta}<{\alpha}\}\subs \br {\omega};{\omega};$,
for $|{\alpha}|\le 2^{\omega}$.
Applying  the inductive hypothesis for each 
${\beta}<{\alpha}$ we also fix a locally compact scattered space 
$X_{\beta}$ of height ${\beta}$ such that $\lev
0;X_{\beta};=A_{\beta}$
and $X_{\beta}\cap X_{\gamma}=A_{\beta}\cap A_{\gamma}$ 
for $\{{\beta},{\gamma}\}\in \br {\alpha};2;$.
Now amalgamate the spaces $X_{\beta}$ as follows:
consider the topological space 
$X=\<\cup_{{\beta}<{\alpha}}X_{\beta},{\tau}\>$ where
${\tau}$ is the topology generated by 
$\cup_{{\beta}<{\alpha}}{\tau}_{X_{\beta}}$.
Since $A_{\beta}\cap A_{\gamma}$ is a finite and open subspace of both 
  $X_{\beta}$ and $X_{\gamma}$ 
it follows that each $X_{\beta}$ is an 
open subspace of $X$. Consequently, $X$ is LCS 
with countably many isolated points, and 
$\htt(X)=\sup_{{\beta}<{\alpha}}\htt X_{\beta}={\alpha}$.
\end{proof}

\begin{corollary}\label{tm:zfc}
There is a locally compact, scattered space of height 
$\oot$ and having $\oo$ isolated points.
\end{corollary}

\begin{proof}
If $2^{\omega}={\omega}_1$, then theorem \ref{tm:scch} gives such a space.

If $2^{\omega}>{\omega}_1$  then $(2^{\omega})^+\ge \ooh$ and so
 according to theorem \ref{tm:notch}
for each ${\alpha}<\ooh$ there is
locally compact, scattered space of height 
${\alpha}$ and countably many isolated points. 
\end{proof}

\end{document}